\documentclass[A4,twoside,11pt,leqno]{article}
\usepackage{amssymb}
\usepackage{amsthm}
\usepackage[tbtags]{amsmath}
\usepackage{doc}
\usepackage{latexsym}
\usepackage{amscd}
\usepackage[frame,cmtip,arrow,matrix,line,graph,curve]{xy}
\usepackage{graphics}
\usepackage{epsfig}
\usepackage{eucal}
\usepackage{mathrsfs}
\font\headd=cmr8

\begin{document}
\thispagestyle{plain}
 \markboth{}{}
\small{\addtocounter{page}{0} \pagestyle{plain}
\noindent{\scriptsize KYUNGPOOK Math. J. 00(0000), 000-000}\\
\noindent{\scriptsize http://dx.doi.org/00.0000/KMJ.0000.00.0.000}\\
\noindent {\scriptsize pISSN 1225-6951 \,\,\,\,\,\, eISSN 0454-8124}\\
\noindent {\scriptsize $\copyright$ Kyungpook Mathematical Journal}
\vspace{0.2in}\\
%%%%%%%%%%%%%%%%%%%%%%%%%%%%%%%%%%%%%%%%%%%%%%%%%%%%%%%%%%%%%%%%%%%%%%%%%%%%%%%%%%
%%%%%%%%%%%%%%%%%%%%%%%%%%%%%%%%%%%  Title in First page  %%%%%%%%%%%%%%%%%%%%%%%%
%%%%%%%%%%%%%%%%%%%%%%%%%%%%%%%%%%%%%%%%%%%%%%%%%%%%%%%%%%%%%%%%%%%%%%%%%%%%%%%%%%
\noindent{\large\bf Equivalence of cyclic $p$-squared actions on handlebodies }
\footnote{{}\\ \\[-0.7cm]
Received March 00, 2014; revised May 00, 2014; accepted November 00, 2014.\\
%%%%%%%%%%%%%%%%%%%%%%%%%%%%%%%%%%%%%%%%%%%%%%%%%%%%%%%%%%%%%%%%%%%%%%%%%%%%%%%%%%
%%%%%%%%%%%%%%%%%%%% 2010 Mathematics Subject Classification %%%%%%%%%%%%%%%%%%%%%
%%%%%%%%%%%%%%%%%%%%%%%%%%%%%%%%%%%%%%%%%%%%%%%%%%%%%%%%%%%%%%%%%%%%%%%%%%%%%%%%%%
%%%%%%%%%%%%%%%%%%%%%%%%%%%%%%%%%%%%%%%%%%%%%%%%%%%%%%%%%%%%%%%%%%%%%%%%%%%%%%%%%%
%%%%%%%%%%%%%%%%%%%%%%%%%%%%%%% Key words and phrases %%%%%%%%%%%%%%%%%%%%%%%%%%%%
%%%%%%%%%%%%%%%%%%%%%%%%%%%%%%%%%%%%%%%%%%%%%%%%%%%%%%%%%%%%%%%%%%%%%%%%%%%%%%%%%%
}
\vspace{0.15in}\\
%%%%%%%%%%%%%%%%%%%%%%%%%%%%%%%%%%%%%%%%%%%%%%%%%%%%%%%%%%%%%%%%%%%%%%%%%%%%%%%%%%
%%%%%%%%%%%%%%%%%%%%%%%%%% Authors & Address & E-mail %%%%%%%%%%%%%%%%%%%%%%%%%%%%
%%%%%%%%%%%%%%%%%%%%%%%%%%%%%%%%%%%%%%%%%%%%%%%%%%%%%%%%%%%%%%%%%%%%%%%%%%%%%%%%%%
\noindent{\sc Jesse Prince-Lubawy}
\newline
{\it Department of Mathematics, University of North Alabama, Florence, Alabama, USA\\
e-mail} : {\verb|jprincelubawy@una.edu|}
\vspace{0.15in}\\
%%%%%%%%%%%%%%%%%%%%%%%%%%%%%%%%%%%%%%%%%%%%%%%%%%%%%%%%%%%%%%%%%%%%%%%%%%%%%%%%%%
%%%%%%%%%%%%%%%%%%%%%%%%%% Abstract Abstract Abstract %%%%%%%%%%%%%%%%%%%%%%%%%%%%
%%%%%%%%%%%%%%%%%%%%%%%%%%%%%%%%%%%%%%%%%%%%%%%%%%%%%%%%%%%%%%%%%%%%%%%%%%%%%%%%%%
{\footnotesize {\sc Abstract.} In this paper we consider all orientation-preserving $\mathbb{Z}_{p^2}$-actions on 3-dimensional handlebodies $V_g$ of genus $g>0$ for $p$ an odd prime. To do so, we examine particular graphs of groups $(\Gamma($v$),\mathbf{G(v)})$ in canonical form for some 5-tuple v $=(r,s,t,m,n)$ with $r+s+t+m>0$. These graphs of groups correspond to the handlebody orbifolds $V(\Gamma($v$),{\mathbf{G(v)}})$ that are homeomorphic to the quotient spaces $V_g/\mathbb{Z}_{p^2}$ of genus less than or equal to $g$. This algebraic characterization is used to enumerate the total number of $\mathbb{Z}_{p^2}$-actions on such handlebodies, up to equivalence.
}
\vspace{0.2in}\\
\pagestyle{myheadings}
 \markboth{\headd Jesse Prince-Lubawy$~~~~~~~~~~~~~~~~~~~~~~~~~~~~~~~~~~~~~~~~~~~~~\,$}
 {\headd $~~~~~~~~~~~~~~~~~~~~~~~~~~~~~~~~~~~~~$Cyclic $p$-squared actions}\\
%
%
%
%
%%%%%%%%%%%%%%%%%%%%%%%%%%%%%%%%%%%%%%%%%%%%%%%%%%%%%%%%%%%%%%%%%%%%%%%%%%%%%%%%%%
%%%%%%%%%%%%%%%%%%%%%%%%%% Contents of Section 1 %%%%%%%%%%%%%%%%%%%%%%%%%%%%%%%%%
%%%%%%%%%%%%%%%%%%%%%%%%%%%%%%%%%%%%%%%%%%%%%%%%%%%%%%%%%%%%%%%%%%%%%%%%%%%%%%%%%%
\noindent{\bf 1. Introduction}
\setcounter{equation}{0}
\renewcommand{\theequation}{1.\arabic{equation}}
\vspace{0.1in}\\
\indent
A {$\mathbf{G}$-action} on  a handlebody $V_g$, of genus $g>0$, is a group monomorphism $\phi:{\mathbf{G}}\longrightarrow{}$Homeo$^+(V_g)$, where Homeo$^+(V_g)$ denotes the group of orientation-preserving homeomorphisms of $V_g$. Two actions $\phi_1$ and $\phi_2$ on $V_g$ are said to be equivalent if and only if there exists an orientation-preserving homeomorphism $h$ of $V_g$ such that $\phi_2(x)=h\circ\phi_1(x)\circ{h^{-1}}$ for all $x\in\mathbf{G}$. From \cite{MMZ}, the action of any finite group $\mathbf{G}$ on $V_g$ corresponds to a collection of graphs of groups.  We may assume these particular graphs of groups are in canonical form and satisfy a set of normalized conditions, which can be found in \cite{KM}.  

Let v\ $=(r,s,t,m,n)$ be an ordered $5$-tuple of nonnegative integers. The graph of groups $(\Gamma($v$),\mathbf{G(v)})$ in canonical form (see \cite{JPL} for the case $p=2$ and $p^2=4$) determines a handlebody orbifold $V(\Gamma($v$),\mathbf{G(v)})$. The orbifold $V(\Gamma($v$),{\mathbf{G(v)}})$ is constructed in a similar manner as described in \cite{KM}. Note that the quotient of any $\mathbb{Z}_{p^2}$-action on $V_g$ is an orbifold of this type, up to homeomorphism. 

An explicit combinatorial enumeration of orientation-preserving $\mathbb{Z}_{4}$-actions on $V_g$, up to equivalence, is given in \cite{JPL}. In this work we will be interested in examining the orientation-preserving geometric group actions on $V_g$ for the group $\mathbb{Z}_{p^2}$ for $p$ an odd prime.

{\bf{Acknowledgement.}} I would like to thank John Kalliongis for his thoughtful advice on polishing my paper and Daniel Branscomb for his insight on pair counting found in Section 3 of the paper.
\vspace{0.2in}\\
%
%
%
%
%
%
%
%
%
%
%
%
%
%%%%%%%%%%%%%%%%%%%%%%%%%%%%%%%%%%%%%%%%%%%%%%%%%%%%%%%%%%%%%%%%%%%%%%%%%%%%%%%%%%
%%%%%%%%%%%%%%%%%%%%%%%%%% Contents of Section 2 %%%%%%%%%%%%%%%%%%%%%%%%%%%%%%%%%
%%%%%%%%%%%%%%%%%%%%%%%%%%%%%%%%%%%%%%%%%%%%%%%%%%%%%%%%%%%%%%%%%%%%%%%%%%%%%%%%%%
\setcounter{equation}{0}
\renewcommand{\theequation}{2.\arabic{equation}}
\noindent{\bf 2. Combinatorial Argument}
\vspace{0.1in}\\
\indent
Let $k\in\mathbb{N}$. Define $S\left(k\right)=\{y_1,y_2,y_3,\ldots,y_k\}$ to be the set such that $y_i$, $i\in\mathbb{N}$ and $1\leq{i}\leq{k}$, is an element of a finite group $G$. Let $y_{i_j}\in S\left(k\right)$, where $1\leq i_j\leq k$, $j\in\mathbb{N}$ and $1\leq{j}\leq{k}$.  Define $C\left(j\right)=\{\left(y_{i_1},y_{i_2},y_{i_3},\ldots,y_{i_j}\right)|\ i_1\leq i_2\leq i_3\leq\cdots\leq i_j\}\subseteq \left(S(k)\right)^j$ to be the set of ordered $j$-tuples. Define $C\left(j,l\right)=\{\left(y_{i_1},y_{i_2},y_{i_3},\ldots,y_{i_j}\right)|\ i_1=k-\left(l-1\right)\}\subseteq C\left(j\right)$ to be the set of ordered $j$-tuples with a fixed first coordinate. We shall use the following lemma to discuss the order of the set $C\left(j\right)$.\vspace{0.1in}\\
\noindent
{\bf Lemma 2.1.}
{\it 
$C\left(j,l\right)\cap C\left(j',l\right)=\emptyset$ if and only if $j\neq j'$.
}\vspace{0.1in}\\
\indent
Now define $\mathscr{C}\left(j\right)=|C\left(j\right)|$ and $\mathscr{C}\left(j,l\right)=|C\left(j,l\right)|$. Then $\mathscr{C}\left(j\right)=\sum_{l=1}^k\mathscr{C}\left(j,l\right)$. Furthermore, the following lemma holds.\vspace{0.1in}\\
\noindent
{\bf Lemma 2.2.}
{\it 
$\mathscr{C}\left(j+1,l\right)=\sum_{u=1}^{l}\mathscr{C}\left(j,u\right)$.
}\vspace{0.1in}\\
\noindent
{\bf Theorem 2.3.}
{\it 
$\mathscr{C}\left(j\right)=A_{k_j}$ where $A_{k_1}=k$ if $j=1$, $A_{k_2}=\frac{k(k+1)}{2}$ if $j=2$, and $A_{k_j}=\sum_{i=0}^{k-1}\left[{j-3+i\choose j-3}\sum_{u=1}^{k-i}u\right]$ if $j\geq3$.
}\vspace{0.05in}\\
\noindent
{\it Proof.}
Let $j=1$. Then $\mathscr{C}(1)=\sum_{l=1}^k\mathscr{C}\left(1,l\right)=\sum_{l=1}^k1=k$. Now let $j=2$. Then $\mathscr{C}\left(2,l\right)=l$ so that $\mathscr{C}(2)=\sum_{l=1}^k\mathscr{C}\left(2,l\right)=\sum_{l=1}^kl=\frac{k(k+1)}{2}$. Finally let $j\geq3$. Then $\mathscr{C}\left(j\right)=\sum_{i=0}^{k-1}\left[{j-3+i\choose j-3}\sum_{u=1}^{k-i}\right]$. To see this we will use induction on $j$. For the base case, let $j=3$. Then
\begin{align}
\mathscr{C}\left(3\right)&=\sum_{l=1}^k\mathscr{C}\left(3,l\right)=\sum_{l=1}^k\sum_{u=1}^l\mathscr{C}\left(2,u\right)=\sum_{l=1}^k\sum_{u=1}^lu=\sum_{u=1}^1u+\sum_{u=1}^2u+\cdots+\sum_{u=1}^ku \notag\\
&=\sum_{u=1}^ku+\cdots+\sum_{u=1}^2u+\sum_{u=1}^1u=\sum_{i=0}^{k-1}\left[{i\choose 0}\sum_{u=1}^{k-i}u\right]. \notag
\end{align}
\noindent
For the inductive step, assume $\mathscr{C}\left(j\right)=\sum_{i=0}^{k-1}\left[{j-3+i\choose j-3}\sum_{u=1}^{k-i}u\right]$. Then 

\begin{align}
\mathscr{C}\left(j+1\right)&=\sum_{l=1}^k\mathscr{C}\left(j+1,l\right)=\mathscr{C}\left(j+1,1\right)+\mathscr{C}\left(j+1,2\right)+\cdots+\mathscr{C}\left(j+1,k\right) \notag\\
&=\sum_{q=1}^1\mathscr{C}\left(j,q\right)+\sum_{q=1}^2\mathscr{C}\left(j,q\right)+\cdots+\sum_{q=1}^k\mathscr{C}\left(j,q\right)=\sum_{i=0}^0\left[{j-3+i\choose j-3}\sum_{u=1}^{1-i}u\right] \notag\\
&+\sum_{i=0}^1\left[{j-3+i\choose j-3}\sum_{u=1}^{2-i}u\right]+\cdots+\sum_{i=0}^{k-1}\left[{j-3+i\choose j-3}\sum_{u=1}^{k-i}u\right] \notag\\
&=\left[{j-3\choose j-3}+{j-3+1\choose j-3}+{j-3+2\choose j-3}+\cdots+{j-3+\left(k-1\right)\choose j-3}\right]\sum_{u=1}^1u \notag\\
&+\left[{j-3\choose j-3}+{j-3+1\choose j-3}+\cdots+{j-3+\left(k-2\right)\choose j-3}\right]\sum_{u=1}^2+\cdots \notag\\
&=\left[{j-2+\left(k-1\right)\choose j-2}\right]\sum_{u=1}^1u. \notag
\end{align}
\hfill$\Box$

\noindent
Note that we will define $A_{k_0}=1$ for $j=0$.
\medskip

%
%
%
%
%
%
%
%
%
%
%
%
%%%%%%%%%%%%%%%%%%%%%%%%%%%%%%%%%%%%%%%%%%%%%%%%%%%%%%%%%%%%%%%%%%%%%%%%%%%%%%%%%%
%%%%%%%%%%%%%%%%%%%%%%%%%% Contents of Section 3 %%%%%%%%%%%%%%%%%%%%%%%%%%%%%%%%%
%%%%%%%%%%%%%%%%%%%%%%%%%%%%%%%%%%%%%%%%%%%%%%%%%%%%%%%%%%%%%%%%%%%%%%%%%%%%%%%%%%
\setcounter{equation}{0}
\renewcommand{\theequation}{3.\arabic{equation}}
\noindent{\bf 3. The 5-tuple v$=(r,s,t,m,n)$ with $s+t>0$}
\vspace{0.2in}\\
\indent
The orbifold fundamental group of $V(\Gamma($v$,{\mathbf{G(v)}}))$ is an extension of $\pi_1(V_g)$ by the group ${\mathbf{G}}=\mathbb{Z}_{p^2}$. We may view the fundamental group as a free product $G_1*G_2*G_3*\cdots*G_{r+s+t+m+n}$, where $G_i$ is isomorphic to either $\mathbb{Z}$, $\mathbb{Z}_{p^2}\times\mathbb{Z}$, $\mathbb{Z}_{p^2}$, $\mathbb{Z}_p\times\mathbb{Z}$, or $\mathbb{Z}_p$. We establish notation similar to \cite{JPL} for the generators of $\pi_1^{orb}(V(\Gamma($v$,{\mathbf{G(v)}})))$. Furthermore, from \cite{JPL}, we may define an equivalence relation and, using techniques found in \cite{KM}, prove the following.\vspace{0.1in}\\
{\bf Lemma 3.1.}
{\it Let ${\mathbf{v}}=(r,s,t,m,n)$. The set $\mathscr{E}(\mathbb{Z}_{p^2},V_g,V(\Gamma({\bf{v}}),{\bf{G(v)}}))$ is in one-to-one correspondence with the set $\Delta(\mathbb{Z}_{p^2},V_g,V(\Gamma({\bf{v}}),{\bf{G(v)}}))$ for every $g>0$.
}\vspace{0.1in}\\
\indent
To prove the three main theorems in this paper, we count the number of elements in the delta set and use the one-to-one correspondence given in Lemma 3.1 to give the total count for the set $\mathscr{E}(\mathbb{Z}_{p^2},V_g,V(\Gamma({\bf{v}}),{\bf{G(v)}}))$. We resort to the following lemma to help count the number of elements in the delta set. The proof is an adaptation from \cite{KM}.\vspace{0.1in}\\
{\bf Lemma 3.2.}
{\it If $\alpha\in$ $Aut(\pi_1^{orb}(V(\Gamma({\bf{v}}),\bf{G(v)})))$, then $\alpha=h_*$ for some orientation-preserving homeomorphism $h:V(\Gamma({\bf{v}}),{\bf{G(v)}})\longrightarrow{V(\Gamma({\bf{v}}),{\bf{G(v)}})}$ if and only if 
\begin{align}
\alpha(b_j)&=x_jb_{\sigma(j)}^{\varepsilon_j}x_j^{-1},\notag \\
\alpha(c_j)&=x_jb_{\sigma(j)}^{{v_j}}c_{\sigma(j)}^{\varepsilon_j}x_j^{-1},\notag \\
\alpha(d_k)&=y_kd_{\tau(k)}^{\delta_k}y_k^{-1},\notag \\
\alpha(e_l)&=u_le_{\gamma(l)}^{\varepsilon_l'}u_l^{-1},\notag \\
\alpha(f_l)&=u_le_{\gamma(l)}^{{w_l}}f_{\gamma(l)}^{\varepsilon_l'}u_l^{-1}, and\notag \\
\alpha(g_q)&=z_qg_{\xi(q)}^{\delta_q'}z_q^{-1},\notag
\end{align}
for some $x_j, y_k, u_l, z_q\in{\pi_1^{orb}(V(\Gamma({\bf{v}}),\bf{G(v)}))}$; $\sigma\in{\sum_s}$, $\tau\in{\sum_t}$, $\gamma\in{\sum_m}$, $\xi\in{\sum_n}$; $\varepsilon_j, \delta_k, \varepsilon_l', \delta_q'\in{\{+1,-1\}}$; and $0\leq{v_j}<p^2$, $0\leq{w_l}<p$.\\
Note that $\Sigma_l$ is the permutation group on l letters.
}\vspace{0.1in}\\
\indent
Note that from \cite{FR}, a generating set for the automorphisms of $\pi_1^{orb}(V(\Gamma($v$,{\mathbf{G(v)}})))$ is the set of mappings $\{\rho_{ji}(x), \lambda_{ji}(x), \mu_{ji}(x), \omega_{ij}, \sigma_i, \phi_i\}$ whose definitions may be found in \cite{FR}. The first five maps are realizable. The realizable $\phi_i$'s are of the form found in Lemma 3.2 and will be used in the remaining arguments of this paper.\vspace{0.1in}\\
{\bf Lemma 3.3.}
{\it Let ${\mathbf{v}}=(r,s,t,m,n)$ for $s+t>0$ and let $\lambda:\pi_1^{orb}(V(\Gamma($v$,{\mathbf{G(v)}})))\longrightarrow\mathbb{Z}_{p^2}$ be a finite injective epimorphism. There exists a finite injective epimorphism $\lambda':\pi_1^{orb}(V(\Gamma($v$,{\mathbf{G(v)}})))\longrightarrow\mathbb{Z}_{p^2}$ equivalent to $\lambda$ such that the following hold:
\begin{enumerate}
\item $\lambda'(a_1)=\cdots=\lambda'(a_r)=0$.
\item $1\leq\lambda'(b_1)=x_1\leq\lambda'(b_2)=x_2\leq\cdots\leq\lambda'(b_s)=x_s\leq\frac{p^2-1}{2}$ and gcd$(x_i,p^2)=1$ for $1\leq i\leq s$.
\item $\lambda'(c_1)=\cdots=\lambda'(c_s)=0$.
\item $1\leq\lambda'(d_1)=y_1\leq\lambda'(d_2)=y_2\leq\cdots\leq\lambda'(d_t)=y_t\leq\frac{p^2-1}{2}$ and gcd$(y_j,p^2)=1$ for $1\leq j\leq t$.
\item $p\leq\lambda'(e_l)=u_l\leq(\frac{p-1}{2})p$ and gcd$(u_l,p^2)=p$ for all $1\leq l\leq m$.
\item $0\leq\lambda'(f_l)\leq p-1$ for all $1\leq l\leq m$.
\item $p\leq\lambda'(g_1)=z_1\leq\lambda'(g_2)=z_2\leq\cdots\leq\lambda'(g_n)=z_n\leq(\frac{p-1}{2})p$ and gcd$(z_k,p^2)=p$ for $1\leq k\leq n$.
\end{enumerate}
}
\noindent
{\it Proof.}
Let $\lambda:\pi_1^{orb}(V(\Gamma($v$,{\mathbf{G(v)}})))\longrightarrow\mathbb{Z}_{p^2}$ be a finite injective epimorphism. Without loss of generality, assume that $s>0$. Then there exists an element $k_i$ such that $k_i\lambda(b_{r+1})=\lambda(a_i)$. Note that if $s=0$, choose $d_{r+1}$. Property (1) follows by composing $\lambda$ with the realizable automorphism $\prod\lambda_{(r+1)i}((\lambda(b_{r+1}))^{k_i})$. (Note that this is a handle slide). Similarly, there exists an element $\ell_i$ such that $\ell_i\lambda(b_i)=\lambda(c_i)$. Property (3) follows by composing $\lambda$ with the realizable automorphism $\prod\phi_i$, where $\phi_i$ sends the generator $c_i$ to the element $b_i^{-\ell_i}c_i$ and fixes $b_i$. (Note that this is a Dehn twist). Now cut the set $\mathbb{Z}_{p^2}-\{0\}$ in half to get the two sets $\{1,2,\ldots,\frac{p^2-1}{2}\}$ and $\{\frac{p^2+1}{2},\ldots,p^2-1\}$. Notice that each element in the first set is the inverse of an element in the second. Property (2) follows by composing $\lambda$ with the realizable automorphism $\prod\phi_i$, where $\phi_i$ sends the generator $b_i$ to $b_i^{-1}$ provided $\lambda(b_i)\in\{\frac{p^2+1}{2},\ldots,p^2-1\}$. (Note that this is a spin). We may then compose with the realizable automorphism $\prod\omega_{ij}$, which interchanges handles if necessary. A similar argument shows Property (4) and (7). Property (5) follows from a spin, if necessary. Now since $|\lambda(e_i)|=p$, we may apply a Dehn twist, if necessary, to ensure that $0\leq\lambda(f_i)\leq p-1$.\vspace{0.05in}\\

\hfill$\Box$

\vspace{0.1in}
\noindent
{\bf Theorem 3.4.}
{\it Let {\bf v}$=(r,s,t,m,n)$ with $s+t>0$. If $\mathbb{Z}_{p^2}$ acts on $V_g$, then $g-1=p^2\left(r+s+m-1\right)+\left(p^2-1\right)t+\left(p^2-p\right)n$. The number of equivalence classes of $\mathbb{Z}_{p^2}$-actions on $V_g$ such that $V_g/\mathbb{Z}_{p^2}=V\left(\Gamma({\bf v}),G({\bf v})\right)$ is the product $A_{(\frac{p(p-1)}{2})_s}\cdot A_{(\frac{p(p-1)}{2})_t}\cdot A_{(\frac{p(p-1)}{2})_m}\cdot A_{(\frac{p-1}{2})_n}$.
}\vspace{0.1in}\\
\noindent
{\it Proof.}
Applying Lemma 3.1, and noting that the orbifold fundamental group is a free product, we see that the count of the delta set is the product of the count of distinct mappings $\lambda$ satisfying the result of Lemma 3.3. Due to this, we only need to consider the generators $b_i$, $d_i$, $e_i$, $f_i$, and $g_i$. We will first count the number of distinct mappings restricted to the generator $b_i$. To do so, note that there are $\frac{p(p-1)}{2}$ generators of $\mathbb{Z}_{p^2}$ in the set $\{1,2,\ldots,\frac{p^2-1}{2}\}$. Since $\lambda(b_i)$ are ordered, we may think of them as $s$-tuples in the set $C(\frac{p(p-1)}{2})$. Applying Theorem 2.3, we see that there are $A_{(\frac{p(p-1)}{2})_s}$ distinct $\lambda$'s. A similar argument works for the generators $d_i$ and $g_i$. (Note that there are only $\frac{p-1}{2}$ multiples of $p$ in the set $\{1,2,\ldots,\frac{p^2-1}{2}\}$). Now to count the number of distinct mappings for the generators $e_i$ and $f_i$, we will note that $(e_i,f_i)\in\mathbb{Z}_{p^2}\times\mathbb{Z}_{p^2}$. When $m=1$, there are $\frac{p(p-1)}{2}$ distinct order pairs that satisfy Properties (5) and (6) from Lemma 3.3. We may relabel these ordered pairs by $y_i$ and create an $m$-tuple in the set $C(\frac{p(p-1)}{2})$. Note that there exists a realizable action when doing this. Again, applying Theorem 2.3, we see that there are $A_{(\frac{p(p-1)}{2})_m}$ distinct $\lambda$'s.
\hfill$\Box$
\medskip

%
%
%
%
%
%
%
%
%
%
%
%
%%%%%%%%%%%%%%%%%%%%%%%%%%%%%%%%%%%%%%%%%%%%%%%%%%%%%%%%%%%%%%%%%%%%%%%%%%%%%%%%%%
%%%%%%%%%%%%%%%%%%%%%%%%%% Contents of Section 4 %%%%%%%%%%%%%%%%%%%%%%%%%%%%%%%%%
%%%%%%%%%%%%%%%%%%%%%%%%%%%%%%%%%%%%%%%%%%%%%%%%%%%%%%%%%%%%%%%%%%%%%%%%%%%%%%%%%%
\setcounter{equation}{0}
\renewcommand{\theequation}{4.\arabic{equation}}
\noindent{\bf 4. The 5-tuple v$=(r,0,0,m,n)$ with $r>0$}
\vspace{0.2in}\\
\indent
We will now consider the 5-tuple v$=(r,0,0,m,n)$ with $s+t=0$ and $r>0$. For the following lemma, we will need to account for two cases: (1) there exists at least one generator $f_i$ that is mapped to a generator of $\mathbb{Z}_{p^2}$ and (2) otherwise. Note that these maps are not equivalent. (This is a modification of Lemma 2.3 from \cite{JPL}).\vspace{0.1in}\\
{\bf Lemma 4.1.}
{\it Let ${\mathbf{v}}=(r,0,0,m,n)$ for $r>0$ and let $\lambda:\pi_1^{orb}(V(\Gamma($v$,{\mathbf{G(v)}})))\longrightarrow\mathbb{Z}_{p^2}$ be a finite injective epimorphism. There exists a finite injective epimorphism $\lambda':\pi_1^{orb}(V(\Gamma($v$,{\mathbf{G(v)}})))\longrightarrow\mathbb{Z}_{p^2}$ equivalent to $\lambda$ such that the following hold:
\begin{enumerate}
\item $\lambda'(a_1)=\cdots=\lambda'(a_r)=0$.
\item $p\leq\lambda'(e_l)=u_l\leq(\frac{p-1}{2})p$ and gcd$(u_l,p^2)=p$ for all $1\leq l\leq m$.
\item $1\leq\lambda'(f_1)=y_1\leq p-1$.
\item $0\leq\lambda'(f_j)\leq p-1$ for all $2\leq j\leq m$.
\item $p\leq\lambda'(g_1)=z_1\leq\lambda'(g_2)=z_2\leq\cdots\leq\lambda'(g_n)=z_n\leq(\frac{p-1}{2})p$ and gcd$(z_k,p^2)=p$ for $1\leq k\leq n$.
\end{enumerate}
OR
\begin{enumerate}
\item $1\leq\lambda'(a_1)=x_1\leq\frac{p(p-1)}{2}$ and gcd$(x_1,p^2)=1$.
\item $\lambda'(a_2)=\cdots=\lambda'(a_r)=0$.
\item $p\leq\lambda'(e_l)=u_l\leq(\frac{p-1}{2})p$ and gcd$(u_l,p^2)=p$ for all $1\leq l\leq m$.
\item $\lambda'(f_1)=\cdots=\lambda'(f_m)=0$.
\item $p\leq\lambda'(g_1)=z_1\leq\lambda'(g_2)=z_2\leq\cdots\leq\lambda'(g_n)=z_n\leq(\frac{p-1}{2})p$ and gcd$(z_k,p^2)=p$ for $1\leq k\leq n$.
\end{enumerate}
}

Note that the first set of properties hold when there exists at least one generator $f_i$ that is mapped to a generator of $\mathbb{Z}_{p^2}$. To show Properties (1)-(5) in both cases, we use similar techniques found in the proof of Lemma 3.3. That is, handle slides, spins, interchanging handles, and Dehn twists. From this we get the following theorem.

\vspace{0.1in}
\noindent
{\bf Theorem 4.2.}
{\it Let {\bf v}$=(r,0,0,m,n)$ with $r>0$. If $\mathbb{Z}_{p^2}$ acts on $V_g$, then $g-1=p^2\left(r+m-1\right)+\left(p^2-p\right)n$. The number of equivalence classes of $\mathbb{Z}_{p^2}$-actions on $V_g$ such that $V_g/\mathbb{Z}_{p^2}=V\left(\Gamma({\bf v}),G({\bf v})\right)$ is the sum $\frac{(p-1)^2}{2}\cdot A_{(\frac{p(p-1)}{2})_{m-1}}\cdot A_{(\frac{p-1}{2})_n}+\frac{p(p-1)}{2}\cdot A_{(\frac{p-1}{2})_m}\cdot A_{(\frac{p-1}{2})_n}$.
}\vspace{0.1in}\\
\noindent
{\it Proof.}
The first portion of the sum follows from the case there exists at least one generator $f_i$ that is mapped to a generator of $\mathbb{Z}_{p^2}$. In this case we can again think of $(e_i,f_i)\in\mathbb{Z}_{p^2}\times\mathbb{Z}_{p^2}$. We may relabel as $y_i$, noting that there are $\frac{(p-1)^2}{2}$ possibilities for $y_1$. The count of the remaining $m-1$ slots follows from Theorem 3.4. (Similarly for counting $g_i$). The second portion of the sum involves counting $a_1$. However, we can see that there are $\frac{p(p-1)}{2}$ possibilities. The remaining values follow from Theorem 3.4.
\hfill$\Box$
\medskip

%
%
%
%
%
%
%
%
%
%
%
%
%%%%%%%%%%%%%%%%%%%%%%%%%%%%%%%%%%%%%%%%%%%%%%%%%%%%%%%%%%%%%%%%%%%%%%%%%%%%%%%%%%
%%%%%%%%%%%%%%%%%%%%%%%%%% Contents of Section 5 %%%%%%%%%%%%%%%%%%%%%%%%%%%%%%%%%
%%%%%%%%%%%%%%%%%%%%%%%%%%%%%%%%%%%%%%%%%%%%%%%%%%%%%%%%%%%%%%%%%%%%%%%%%%%%%%%%%%
\setcounter{equation}{0}
\renewcommand{\theequation}{5.\arabic{equation}}
\noindent{\bf 5. The 5-tuple v$=(0,0,0,m,n)$ with $m>0$}
\vspace{0.2in}\\
\indent
Finally we will consider the 5-tuple v$=(0,0,0,m,n)$ with $r+s+t=0$ and $m>0$. Given a finite injective epimorphism $\lambda$, it is clear that we may obtain an equivalent $\lambda'$ that satisfies Properties (2)-(5) of the first case of Lemma 4.1. The following theorem is a modification of Theorem 4.2.\vspace{0.1in}\\
\noindent
{\bf Theorem 5.1.}
{\it Let {\bf v}$=(0,0,0,m,n)$ with $m>0$. If $\mathbb{Z}_{p^2}$ acts on $V_g$, then $g-1=p^2\left(m-1\right)+\left(p^2-p\right)n$. The number of equivalence classes of $\mathbb{Z}_{p^2}$-actions on $V_g$ such that $V_g/\mathbb{Z}_{p^2}=V\left(\Gamma({\bf v}),G({\bf v})\right)$ is the product $\frac{(p-1)^2}{2}\cdot A_{(\frac{p(p-1)}{2})_{m-1}}\cdot A_{(\frac{p-1}{2})_n}$.
}

\medskip

%
%
%
%
%
%
%
%
%
%
%
%
%%%%%%%%%%%%%%%%%%%%%%%%%%%%%%%%%%%%%%%%%%%%%%%%%%%%%%%%%%%%%%%%%%%%%%%%%%%%%%%%%%
%%%%%%%%%%%%%%%%%%%%%%%%%% Contents of Section 6 %%%%%%%%%%%%%%%%%%%%%%%%%%%%%%%%%
%%%%%%%%%%%%%%%%%%%%%%%%%%%%%%%%%%%%%%%%%%%%%%%%%%%%%%%%%%%%%%%%%%%%%%%%%%%%%%%%%%
\setcounter{equation}{0}
\renewcommand{\theequation}{6.\arabic{equation}}
\noindent{\bf 6. The number of equivalence classes of $\mathbb{Z}_{25}$-actions on $V_{26}$}
\vspace{0.2in}\\
\indent
Fix $p=5$ and $g=26$. Now $g$ must satisfy the genus equation $g=1-|G|\chi(\Gamma,G)$. Therefore we see that $50=24t+20n+25(r+s+m)$. Solving this equation we see that $t=0$, $n=0$, and $r+s+m=2$. This leads to the following ordered 5-tuples: (0,2,0,0,0), (2,0,0,0,0), (0,0,0,2,0), (1,1,0,0,0), (1,0,0,1,0), and (0,1,0,1,0). Using Theorems 3.4, 4.2, and 5.1, the counts for the following ordered 5-tuples are 55, 10, 55, 10, 18, and 100, respectively. Thus the total number of equivalence classes of $\mathbb{Z}_{25}$-actions on $V_{26}$ is 55+10+55+10+18+100=248.

\medskip

\footnotesize{

}

\end{document}